# Graph Theoretical Models of Closed n-Dimensional Manifolds: Digital Models of a Moebius Strip, a Torus, a Projective Plane a Klein Bottle and n-Dimensional Spheres

Alexander V. Evako

"Dianet", Laboratory of Digital Technologies, Moscow, Russia

e-mail: evakoa@mail.ru

**Abstract**

In this paper, we show how to construct graph theoretical models of n-dimensional continuous objects and manifolds. These models retain topological properties of their continuous counterparts. An LCL collection of n-cells in Euclidean space is introduced and investigated. If an LCL collection of n-cells is a cover of a continuous n-dimensional manifold then the intersection graph of this cover is a digital closed n-dimensional manifold with the same topology as its continuous counterpart. As an example, we prove that the digital model of a continuous n-dimensional sphere is a digital n-sphere with at least $2n+2$ points, the digital model of a continuous projective plane is a digital projective plane with at least eleven points, the digital model of a continuous Klein bottle is the digital Klein bottle with at least sixteen points, the digital model of a continuous torus is the digital torus with at least sixteen points and the digital model of a continuous Moebius band is the digital Moebius band with at least twelve points.

## 1. Introduction

In the past decades, non-orientable surfaces such as a Moebius strip, Klein bottle and projective plane have attracted many scientists from different fields. The study was derived from obvious practical and science background. In physics, a considerable interest has emerged in studying lattice models on non-orientable surfaces as new challenging unsolved lattice-statistical problems and as a realization and testing of predictions of the conformal field theory (see, e.g., [14]). Paper [15] studies topological phases on non-orientable spatial surfaces including Möbius strip, Klein bottle, etc. A formulation of gauge theories which applies to non-orientable manifolds is studied in paper [1]. Discrete models of a physical continuous three-dimensional space were investigated in [10]. Since analytic solutions of problems arising in theories involving non-orientable surfaces can be obtained only on simple geometric regions, it is more reasonable to replace non-orientable surfaces with discrete models, that allows to apply computational or numerical methods. Constructing and analyzing n-dimensional digitized images of continuous objects is central to many applications: geoscience, computer graphics, fluid dynamics and medicine [16] are obvious examples. There are a number of different frameworks to represent digital images. Graph-theoretic approach equips a digital image with a graph structure based on the local adjacency relations of points (see e.g. [3]). In paper [4], digital n-surfaces were defined as simple undirected graphs and basic properties of n-surfaces were studied. Properties of digital n-manifolds were investigated in ([2], [5], [6], [9]). Paper [7] shows that the intersection graph of any LCL cover of a plane is a digital 2-dimensional plane. The present paper extents the basic results obtained in [7] to n dimensions. An LCL collection of topological n-disks is a key ingredient for this approach. A digital model for a continuous closed n-dimensional manifold is the intersection graph of an LCL cover of the manifold. It turns out that the digital model is necessarily a digital n-dimensional manifold preserving global and local topological features of a continuous closed n-dimensional manifold.

The material to be presented in Section 2 begins with the investigation of properties of LCL collections of n-cells. In particular, it is proven that if W={$D_0,D_1,…D_s$} is an LCL collection of n-cells then U={$C_1,…C_s$} is an LCL collection of (n-1)-cells, where $C_i=D_0 \cap D_i$, i=1,…s.

Section 3 includes some results related to graphs that are digital n-dimensional manifolds. The examples are given of digital n-dimensional spheres, a digital torus, a digital projective plane, a digital Klein bottle and a digital Moebius band.

Section 4 studies connection between an LCL cover of a closed n-manifold and the intersection graph of this cover. It is shown that the intersection graph of any LCL cover of a continuous closed n-manifold is necessarily a digital closed n-dimensional manifold with the same topological features as its continuous counterpart.

**Computer experiments**

The approach developed in this paper is based on computer experiments described in [12]. Divide Euclidean space $E^n$ into the set of n-cubes with the side length L and vertex coordinates in the set X={ $Lx_1,…Lx_n$: $x_i \in Z$}. Suppose that S is a continuous space in $E^n$. Call the cubical model of S the family M of n-cubes intersecting S, and the digital model of S the intersection graph G of $M_L$.

It was revealed in computer experiments that the length $L_0$ exists such that for any $L_1$ and $L_2<L_0$, spaces S, $M_1$ and $M_2$ are homotopy equivalent, and graphs $G_1$ and $G_2$ can be transformed from one to the other with transformations called contractible in [12].

In other words, we have three mathematical objects:

$$S \rightarrow M(S) \rightarrow G(S)$$

where S is a continuous space, M(S) is a cubical model of S, and G(S) is the intersection graph of M(S). Experiments show that if spaces $S_1$ and $S_2$ are homotopy equivalent (and naturally $M(S_1)$ and $M(S_2)$ are all homotopy equivalent) then graphs $G(S_1)$ and $G(S_2)$ are homotopy equivalent, i.e., can be transformed from one to the other with contractible transformations.

These results yield the following conclusion.
- The digital model (intersection graph) contains topological and perhaps geometrical characteristics of space S. Otherwise, the digital model G of S is a digital counterpart of S.
- Contractible transformations of graphs are a digital analogue of homotopy in algebraic topology.

**Example 1.1**

To illustrate these experiments, consider examples depicted in figure 1. For a curve $S_1$, $M_1$ is a set of cubes intersecting $S_1$, and $G_1$ is the intersection graph of $M_1$. $S_1$ and $M_1$ have the homotopy type of a point, and graph $G_1$ can be transformed with contractible transformations to a one-point graph. $S_2$ is a topological ball, $M_2$ is a set of cubes intersecting $S_2$, and $G_2$ is the intersection graph of $M_2$. As above, $S_2$ and $M_2$ are contractible sets, and graph $G_2$ can be transformed with contractible transformations to a one-point graph.

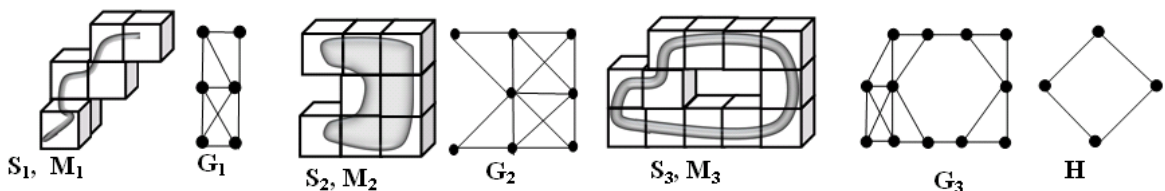

**Figure 1** $S_1$ is a curve, $M_1$ is the cubical model for $S_1$, $G_1$ is the intersection graph of $M_1$. $S_2$ is a topological ball, $M_2$ is the cubical model for $S_2$, $G_2$ is the intersection graph of $M_2$. $S_3$ is a closed curve, $M_3$ is the cubical model for $S_3$, $G_3$ is the intersection graph of $M_3$. H is a digital circle

**Example 1.2**

For a circle $S_3$ shown in figure 1, $M_3$ is a set of cubes intersecting $S_3$, and $G_3$ is the intersection graph of $M_3$. $M_3$ is homotopy equivalent to a unit circle, and $G_3$ can be converted to graph H by contractible transformations. H is a digital 1-dimensional sphere (see the definition below).

## 2. LCL collections of n-cells

In the present paper, we consider closed n-manifolds, which are connected compact n-dimensional manifolds without boundary. That is, a closed n-manifold admits an LCL decomposition. In this section, the primary object we deal with is an LCL collection of topological n-dimensional disks. Notice, that we use intrinsic topology of an object, without reference to an embedding space.

If a set D is homeomorphic to a closed n-dimensional ball on $R^n$, then D is said to be an n-cell or n-ball. A set S is said to be an n-sphere if S is homeomorphic to an n-dimensional sphere $S^n$ on $R^{n+1}$. We denote the interior and the boundary of an n-cell D=Int D∪∂D by IntD and ∂D respectively. ∂D is an (n-1)-sphere. The 0-cell D is a single point for which ∂D=∅. Notice that an n-cell is a contractible space. Remind that collections of sets W={$D_1,D_2,...$} and U={$C_1,C_2,...$} are isomorphic if $D_i \cap D_k \cap ... D_p \neq \emptyset$ if and only if $C_i \cap C_k \cap ... C_p \neq \emptyset$. Or equally: if the intersection graphs G(W) and G(U) of W and U are isomorphic. Facts about n-cells and n-spheres that we will need in this paper are stated below.

**Fact 2.1**

If S is an n-sphere and D is an n-cell contained in S then S-IntD is an n-cell.

**Fact 2.2**

Let $D_1$ and $D_2$ be n-cells such that $D_1 \cap D_2 = \partial D_1 \cap \partial D_2 \neq \emptyset$. The union $D_1 \cup D_2 = B$ is an n-cell iff $D_1 \cap D_2 = D^{n-1}$ is an (n-1)-cell.

Note that the shape and the size of n-cells can be arbitrary. In combinatorics, a collection W={$D_i | i \in I$} is a Helly collection (of order 2) if for every $J \subseteq I$ such that and $D_i \cap D_k \neq \emptyset, i,k \in J$, it follows that $\cap\{D_i | i \in J\} \neq \emptyset$. A locally centered collection of n-cells defined below is a Helly collection.

**Definition 2.1**

A collection W={$D_i | i \in I$} of n-cells is called locally centered (LC) if for every $J \subseteq I$ such that $D_i \cap D_k \neq \emptyset, i,k \in J$, it follows that $\cap\{D_i | i \in J\} \neq \emptyset$.

In paper [8], an LC collection of arbitrary sets on $E^n$ was called continuous.

**Definition 2.2**

Collection of n-cells W is said to be a locally lump (LL) collection if for every non-empty intersection E=$D_1 \cap ... \cap D_k \neq \emptyset$ of k distinct members of W, E is the intersection of their boundaries and is an (n+1-k)-cell, i.e., E=$D_1 \cap ... \cap D_k = \partial D_1 \cap ... \cap \partial D_k = D_{1...k}$ is an (n+1-k)-cell.

**Definition 2.3**

Let W={$D_1,D_2,...D_s$} be a locally centered and locally lump collection of n-cells. Then W is called a locally centered lump (LCL) collection.

If W is an LCL collection of n-cells and $D_1 \cap D_2 \cap ... D_p \neq \emptyset$ then p≤n+1. Figures 2-3 show collections of 1-, 2- and 3-cells. Collections $W_1$, $W_2$ and $W_3$ depicted in figure 2 are non-LCL collections. $W_1$ is not LCL because $D_1 \cap D_2 \not\subseteq \partial D_1 \cap \partial D_2$. $W_2$ is not LCL because $D_1 \cap D_2$ is not a 2-cell. $W_3$ is not LCL because $D_1 \cap D_2 \cap D_3 = \emptyset$. In figure 3, all collections are LCL collections. Notice that an LCL collection is a special case of a complete collection of sets on $E^n$ whose properties were studied in [8]. Paper [7] investigates properties of LCL collections of 1- and 2-cells.

**Remark 2.1**

It follows from definition 2.3 that any subcollection of an LCL collection of n-cells is also an LCL collection.

Structural properties of collections of convex polygons in $R^2$ and polyhedra in $R^3$ were studied in a number of works. Paper [13] presents extensions of the basic results about SN sets of convex polytopes to n dimensions. It was proven that if a set W of convex tiles is strongly normal (SN) then the union N(D) of all Q∈W that intersect D (including D itself) is a contractible set. If the union

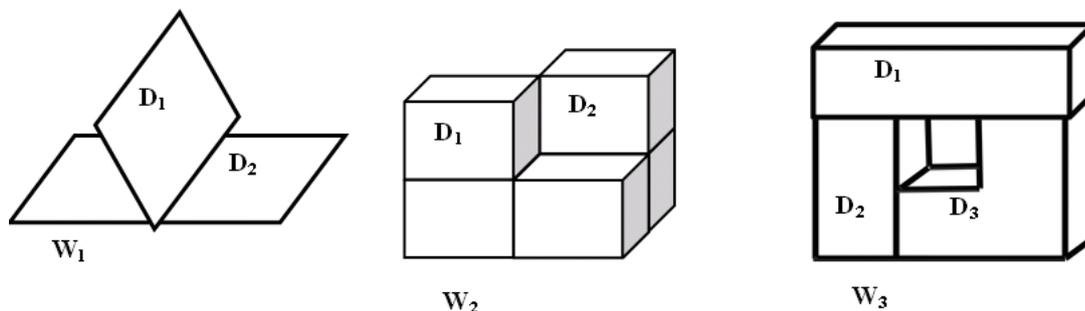

**Figure 2** $W_1$ is not LCL because $D_1 \cap D_2 \not\subseteq \partial D_1 \cap \partial D_2$. $W_2$ is not LCL because $D_1 \cap D_2$ is not a 2-cell. $W_3$ is not LCL because $D_1 \cap D_2 \cap D_3 = \emptyset$.

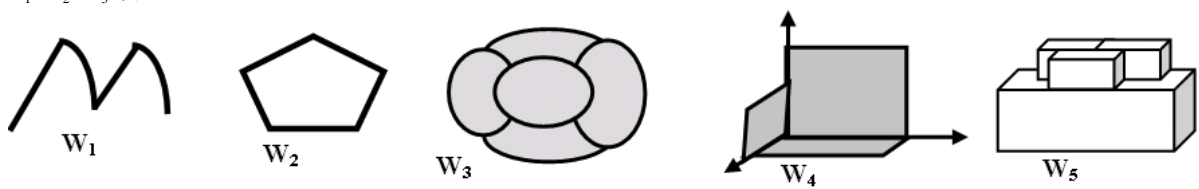

**Figure 3** $W_1 – W_5$ are LCL collections

U(W) of all tiles is contractible then D and the union N(D) are a strong deformation retracts of U(W), i.e. U(W) can be continuously deformed over itself to either D or N(D). If the union $N^*(D)$ of all Q∈W that intersect D (excluding D itself) is contractible then D is called a simple tile, and the union U(W-{D}) of all tiles (excluding D) is a strong deformation retracts of U(W) according to theorem 5.3. It is easy to deduce that an LCL collection of n-cells is SN, but the reverse is not true. All results obtained in [13] are also valid for an LCL collection of n-cells. In particular, if W is an LCL collection of n-cells and the union U(W) is an n-ball, i.e. a contractible set then for any D∈W, U(W) can be continuously deformed over itself to D. This means that W can be transformed into D by deleting simple n-cells.

**Property 2.1**

Let $W=\{D_0, D_1, \ldots D_s\}$ be an LCL collection of n-cells, and the union N(W) of all D is a contractible set (an n-cell). Then there is $D_k$ such that the union M(D) of all D∈W that intersect $D_k$ (excluding $D_k$ itself) is contractible.

This property follows from results of [13], and can be easily verified for any LCL collection with the finite number of n-cells.

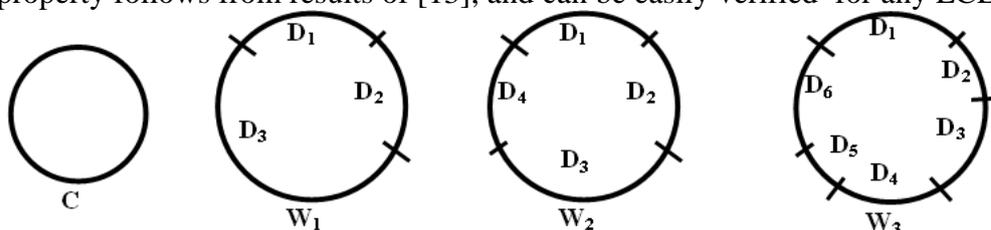

**Figure 4** $W_1$ is a non-LCL cover of a circle C. $W_2$ and $W_3$ are LCL covers of a circle.

The use of LCL collections is rooted in building digital models of continuous objects proposed in this paper. A collection of n-cells intersecting a given n-cell D generates the collection of (n-1)-cells lying on ∂D. The next theorem establishes a general link between the collection of n-cells intersecting D the collection of (n-1)-cells.

**Theorem 2.1**

Let $W=\{D_0,D_1,...D_s\}$ be an LCL collection of n-cells, $V=\{D_1,...D_p\}$ be the collection of n-cells intersecting $D_0$, and $U=\{C_1,...C_p\}$ be the collection of (n-1)-cells which are the intersections $C_i=D_0\cap D_i$, $i=1,...p$. Then U is an LCL collection of (n-1)-cells, and collections V and U are isomorphic.

Proof.

By definition 2.1, any $C_i$ is an (n-1)-cell. To prove that U is an LCL collection, suppose that $C_k\cap C_i\neq\varnothing$, $i,k=1,2,...r$. Since $C_k\cap C_i=D_0\cap D_k\cap D_i$ then $D_0\cap D_i\neq\varnothing$, $D_0\cap D_k\neq\varnothing$, $D_i\cap D_k\neq\varnothing$. Hence,

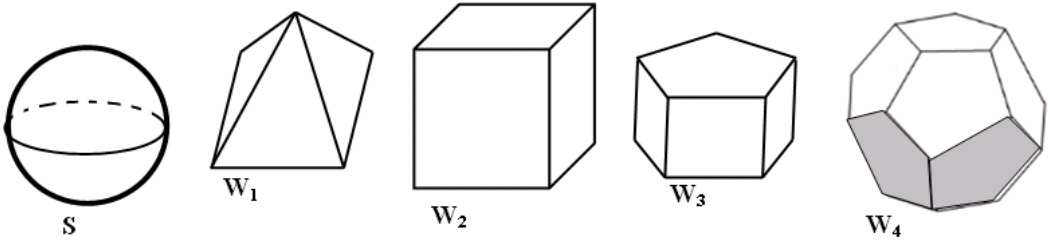

**Figure 5** $W_1$ is a non-LCL cover of a sphere S. The intersection graph $G(W_1)$ is not a digital 2-sphere. $W_2$, $W_3$ and $W_4$ are LCL covers of a sphere. Their intersection graphs $G(W_2)$, $G(W_3)$ and $G(W_4)$ are digital 2-spheres

$C_1\cap C_2\cap...C_r= D_0\cap D_1\cap...D_r=D_{01...r}$ is an (n-r)-cell by definition 2.3. Therefore, U is an LCL collection of (n-1)-cells.

To prove that collections V and U are isomorphic suppose that $C_k\cap C_i=D_0\cap D_k\cap D_{ii}\neq\varnothing$. Then $D_k\cap D_{ii}\neq\varnothing$. Suppose that $C_k\cap C_i=D_0\cap D_k\cap D_i=\varnothing$. Since $D_0\cap D_i\neq\varnothing$, $D_0\cap D_k\neq\varnothing$, then $D_i\cap D_k=\varnothing$. Therefore, collections U and V are isomorphic according to the definition of isomorphic collections of sets. This completes the proof. □

The following result is quite important since allows to reduce the number of elements in an LCL collection.

**Theorem 2.2**

Let $W=\{D_0,D_1,...D_s\}$ be an LCL collection of n-cells intersecting $D_0$, $D_0\cap D_i\neq\varnothing$, for $i=1,...s$. Then $B_k=D_0\cup D_1\cup...D_k$, $k\leq s$, is an n-cell and $V_k=\{B_k,D_{k+1},...D_s\}$ is an LCL collection of n-cells.

Proof.

Let $k=1$. $B_1=D_0\cup D_1$ is an n-cell by fact 2.2. Suppose that $B_1\cap D_i\neq\varnothing$, $D_i\cap D_j\neq\varnothing$, $i,j=2,...p$. We need to show that the intersection $C=B_1\cap D_2\cap...D_p$ is an (n+1-p)-cell.

Since $B_1=D_0\cup D_1$, then $C=(D_0\cup D_1)\cap D_2\cap...D_p= (D_0\cap D_2\cap...D_p)\cup(D_1\cap D_2\cap...D_p)$. Since $D_0\cap D_i\neq\varnothing$, $i=2,...p$, and $D_i\cap D_j\neq\varnothing$, $i,j=2,...p$ then $D_0\cap D_2\cap...D_p=D_{02...p}$ is an (n+1-p)-cell according to definition 2.2.

If $D_1\cap D_2\cap...D_p=\varnothing$ then $C=D_{02...p}$ is an (n+1-p)-cell.

If $D_1\cap D_2\cap...D_p\neq\varnothing$ then $D_1\cap D_2\cap...D_p=D_{12...p}$ is an (n+1-p)-cell. according to definition 2.3. Therefore, $C=D_{02...p}\cup D_{12...p}$. Consider $D_{02...p}\cap D_{12...p}= (D_0\cap D_2\cap...D_p)\cap(D_1\cap D_2\cap...D_p)=D_0\cap D_1\cap D_2\cap...D_p$. Since W is an LCL collection then $D_0\cap D_1\cap D_2\cap...D_p=D_{012...p}$ is an (n-p)-cell according to definition 2.2. Therefore, $C=D_{02...p}\cup D_{12...p}$ is an (n+1-p) cell by fact 2.2. Hence, $V_1=\{B_1,D_2,...D_s\}$ is an LCL collection of n-cells. Sequentially applying the same arguments for $k=2,3,...s$, we obtain that $B_k$ is an n-cell, $V_k=\{B_k,D_{k+1},...D_s\}$, is an LCL collection of n-cells and the union $B_s=D_0\cup D_1\cup...D_s$ is an n-cell. This completes the proof. □

Theorem 2.2 says that the number of elements in an LCL collection can be reduced by merging n-cells provided that the obtained collection is an LCL collection of n-cells. A simple set of n-cells defined below can be replaced by the union of n-cells.

**Definition 2.4**

- Let $W=\{D_0, D_1, \ldots D_s\}$ be an LCL collection of n-cells. A subcollection $W_1=\{D_0, D_1, \ldots D_k\}$ is called simple in W if the union $C_k = D_0 \cup D_1 \ldots \cup D_k$ is an n-cell, and $V=\{C_k, D_{k+1}, \ldots D_s\}$ is an LCL collection n-cells.
- An LCL collection of n-cells is called irreducible (minimal, compressed) if it has no simple

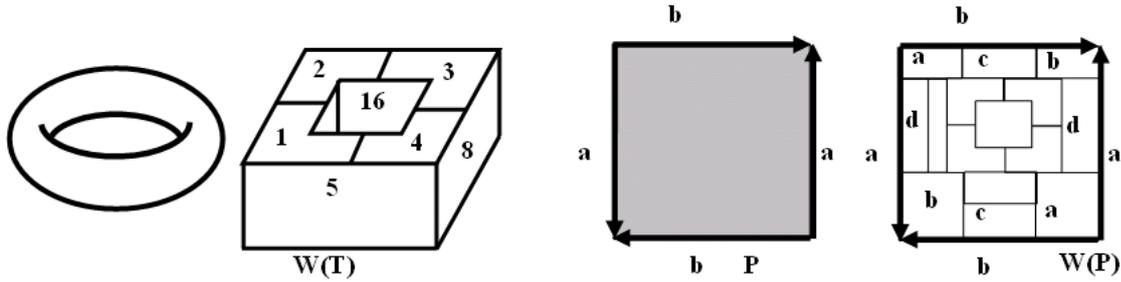

**Figure 6** W(T) is an LCL cover of a torus T. W(P) is an LCL cover of a projective plane P.

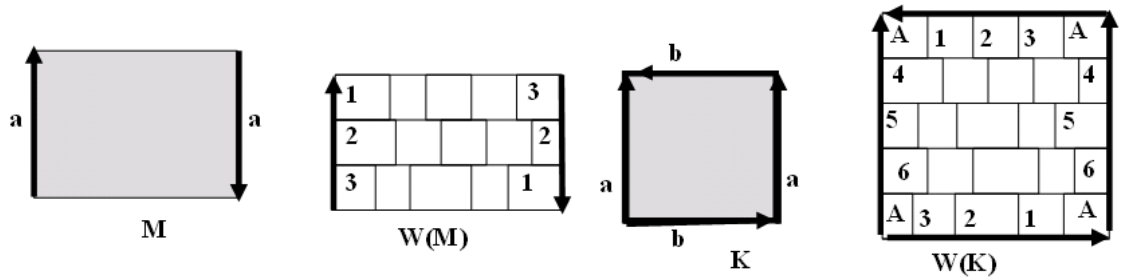

**Figure 7**. M is a continuous Moebius band. W(M) is an LCL cover of M containing twelve elements. K is a continuous Klein bottle. W(K) is an LCL cover of K.

set of n-cells.

In particular, if $W=\{D_0, D_1, \ldots D_s\}$ is an LCL collection of n-cells and the union $U(W) = D_0 \cup D_1 \cup \ldots D_s$ is n-cell (contractible set) then W can be transformed to collection $V=\{U(W)\}$ be merging n-cells contained in simple subcollections.

In figure 4, $W_1$ is not an LCL cover of a circle C, $W_2$ and $W_3$ are LCL covers of C. $W_2$ is an irreducible LCL collection, $W_3$ can be transformed to an irreducible LCL collection by merging $D_3$ and $D_4$, and $D_5$ and $D_6$. Figure 5 depicts LCL covers $W_2$, $W_3$ and $W_4$ of a sphere S. $W_2$ is irreducible, $W_3$ and $W_4$ can be transformed to an irreducible form by merging simple sets of 2-cells. LCL covers of a torus T. and a projective plane P are shown in figure 6, LCL covers of a Moebius band M and a Klein bottle K are shown in figure 7.

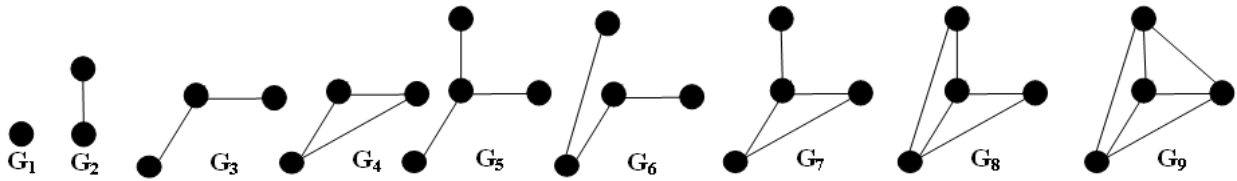

**Figure 8** Contractible graphs with the number of points n<5.

## 3. Digital N-surfaces

There are a considerable amount of literature devoted to the study of different approaches to digital lines, surfaces and spaces in the framework of digital topology. Digital topology studies topological properties of discrete objects which are obtained digitizing continuous objects.

This section includes some results related to digital spaces. Traditionally, a digital image has a graph structure (see [2], [3], [4]). A digital space G is a simple undirected graph $G=(V, W)$ where $V=(v_1, v_2, \ldots v_n, \ldots)$ is a finite or countable set of points, and $W = ((v_p v_q), \ldots)$ is a set of edges. Topological properties of G as a digital space in terms of adjacency, connectedness and dimensionality are completely defined by set W. Let G and v be a graph and a point of G. In ([11], [12]), the subgraph O(v) containing all neighbors of v (without v) is called the rim of point v in G. The subgraph U(v)=v∪O(v) containing O(v) as well as point v is called the ball of point v in G. Let (vu) be an edge of G. The subgraph O(vu)=O(v)∩O(u) is called the rim of (vu).

For two graphs G=(X, U) and H=(Y, W) with disjoint point sets X and Y, their join G⊕H is the graph that contains G, H and edges joining every point in G with every point in H. Contractible graphs are

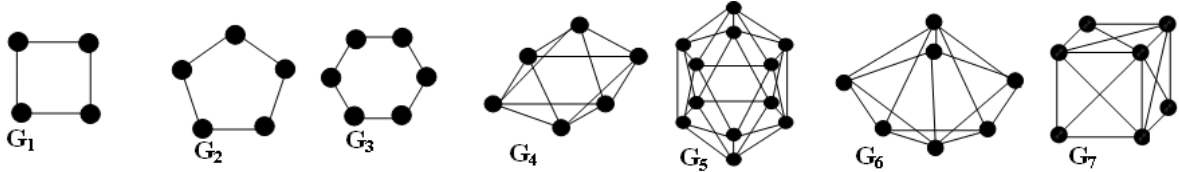

**Figure 9.** $G_1$, $G_2$, $G_3$ are digital one-dimensional spheres. $G_4$, $G_5$, $G_6$ are digital two-dimensional spheres. $G_7$ is a digital three-dimensional sphere

basic elements in this approach.

**Definition 3.1**

- A one-point graph is contractible. If G is a contractible graph and H is a contractible subgraph of G then G can be converted into H by sequential deleting simple points.
- A point v in graph G is simple if the rim O(v) of v is a contractible graph.
- An edge (uv) of a graph G is called simple if the rim O(vu)=O(v)∩O(u) of (uv) is a contractible graph.

By construction, a contractible graph is connected. It follows from definition 1 that a contractible graph can be converted to a point by sequential deleting simple points.

**Definition 3.2**

- Deletions and attachments of simple points and edges are called contractible transformations.
- Graphs G and H are called homotopy equivalent if one of them can be converted to the other one by a sequence of contractible transformations.

Homotopy is an equivalence relation among graphs. Contractible transformations of graphs seem to play the same role in this approach as a homotopy in algebraic topology. In papers [11], [12], it was shown that contractible transformations retain the Euler characteristic and homology groups of a graph. Figure 8 depicts contractible graphs with the number of points n<5. A digital n-manifold is a special case of a digital n-surface defined and investigated in [4].

**Definition 3.3**

- The digital 0-dimensional surface $S^0$(a, b) is a disconnected graph with just two points a and b. For n>0, a digital n-dimensional surface $G^n$ is a nonempty connected graph such that, for each point v of $G^n$, O(v) is a finite digital (n-1)-dimensional surface.
- A connected digital n-dimensional surface $G^n$ is called a digital n-sphere, n>0, if for any point v∈$G^n$, the rim O(v) is an (n-1)-sphere and the space $G^n$-v is a contractible graph (fig.9).
- A digital n-dimensional surface $G^n$ is a digital n-manifold if for each point v of $G^n$, O(v) is a finite digital (n-1)-dimensional sphere (fig.10)
- Digital n-manifolds are called homeomorphic if they are homotopy equivalent. A digital n-manifold M can be converted to a homeomorphic digital n-manifold N with the minimal number of points by contractible transformations.

**Definition 3.4**

- Let M be a digital n-sphere, n>0, and v be a point belonging to M. The space N=M-v is called a digital n-disk with the boundary ∂N=O(v) and the interior IntN=N-∂N.

• Let M be an n-manifold and a point v belong to M. Then the space N=M-v is called an n-manifold with the spherical boundary ∂N=O(v) and the interior IntN=N-∂N.

According to definition 3.3, a digital n-disk is a contractible graph. A digital n-disk is a digital counterpart of a continuous n-dimensional disk in Euclidean $E^n$.

The following results were obtained in [4] and [6].

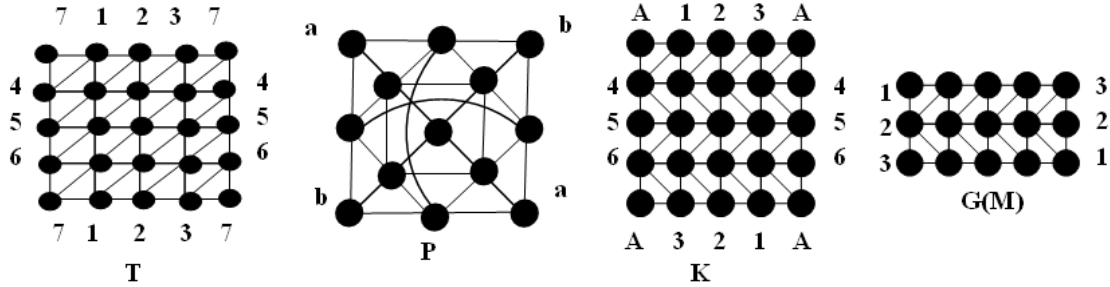

**Figure 10**. T is a digital torus. P is a digital projective plane. K is a digital Klein bottle. M is a digital Moebius band.

**Theorem 3.1**

- The join $S^n_{min}=S^0_1 \oplus S^0_2 \oplus \ldots S^0_{n+1}$ of (n+1) copies of the zero-dimensional sphere $S^0$ is a minimal n-sphere.
- Let M and N be n and m-spheres. Then M⊕N is an (n+m+1)-sphere.
- Any n-sphere M can be converted to the minimal n-sphere $S_{min}$ by contractible transformations.
- Let M be an n-manifold, G and H be contractible subspaces of M and v be a point in M. Then subspaces M-G, M-H and M-v are all homotopy equivalent to each other. The replacement of an edge with a point increases the number of points in a digital n-manifold.

**Definition 3.5**

Let M be an n-manifold, v and u be adjacent points in M and (vu) be the edge in M. Glue a point x to M in such a way that O(x)=v⊕u⊕O(vu), and delete the edge (vu) from the space. This pair of contractible transformations is called the replacement of an edge with a point or R-transformation, R: M→N. The obtained space N is denoted by N=RM=(M∪x)-(vu).

**Theorem 3.2 ([6])**

Let M be an n-manifold and N=RM be a space obtained from M by an R-transformation. Then N is a digital n-manifold homeomorphic to M.

An R-transformation is a digital homeomorphism because it retains the dimension and other local and global topological features of an n-manifold. R-Transformations increase the number of points in a given n-manifold M retaining the global topology (the homotopy type of M) and the local topology (the homotopy type and the dimension of the neighborhood of any point). A close connection between continuous and digital n-manifolds for n=2 was studied in [7].

A digital 0-dimensional surface is a digital 0-dimensional sphere. Figure 9 shows digital 1-, 2- and 3-dimensional spheres. All 2-dimensional spheres are homeomorphic and can be converted into the minimal sphere $S^2_{min}$ by contractible transformations. A digital torus T, a digital projective plane P, a digital Klein bottle K and a digital Moebius band are depicted in figure 10. Notice that paper [10] studies properties of a discrete model of physical continuous three-dimensional space that is a digital three-dimensional Euclidean space.

**4. Connection between an LCL cover of a closed n-manifold and the intersection graph of the LCL cover**

Digital models of continuous objects by using the intersection graphs of covers were studied in [5], [8]. In this section, we are going to show that an LCL cover of a closed n-manifold incorporates topological features directly into the digital model of this manifold. Remind the definition of the intersection graph of a family of sets. Let $W=\{D_1,D_2,....D_n, ...\}$ be a finite or countable family of sets. Then graph $G(W)$ with points $\{x_1,x_2,....x_n,...\}$ is called the intersection graph of W, if points $x_k$ and $x_i$ are adjacent whenever $D_k \cap D_i \neq \emptyset$. In other word, f: $G(W) \to W$ such that $f(x_i)=D_i$ is an isomorphism.

**Definition 4.1**

Let $W=\{D_1,D_2,...D_s\}$ be an LCL collection of n-cells and $M=D_1 \cup D_2 \cup ...D_s$. The intersection graph $G(W)$ of W is called the digital model of M in regard to W.

Notice that by nature, digital models (graphs) of continuous closed n-dimensional manifolds are different from simplicial complexes which are triangulations of the manifolds.

**Theorem 4.1**

Let X be an n-ball, and an LCL collection $W=\{D_1,D_2...D_s\}$ of n-cells be a cover of $X=U(W)=D_1 \cup D_2 \cup ...D_s$. Then the intersection graph $G(W)=\{v_1,v_2,...v_s\}$ of W is contractible.

Proof.
The proof is by induction on the number elements s. For small s, the statement is verified directly. Assume that the statement is valid whenever s<m. Let s=m. Since $X=U(W)$ is contractible then W can be converted into $D_1$ by sequential deleting simple n-cells $D_2,...D_s$ according to remark 2.2. Since $D_2$ is simple in W then the intersection graph $O(v_2)$ of all n-cells that intersect $D_2$ (excluding $D_2$ itself) is contractible by the inductive hypothesis. Therefore, point $v_2$ is simple in $G(W)$. Applying the same arguments we obtain that $G(W)$ is converted into point $v_1$ by sequential deleting simple points. Hence, $G(W)$ is contractible. This completes the proof. □

We will show now that a digital model of an n-dimensional sphere is a digital n-sphere.

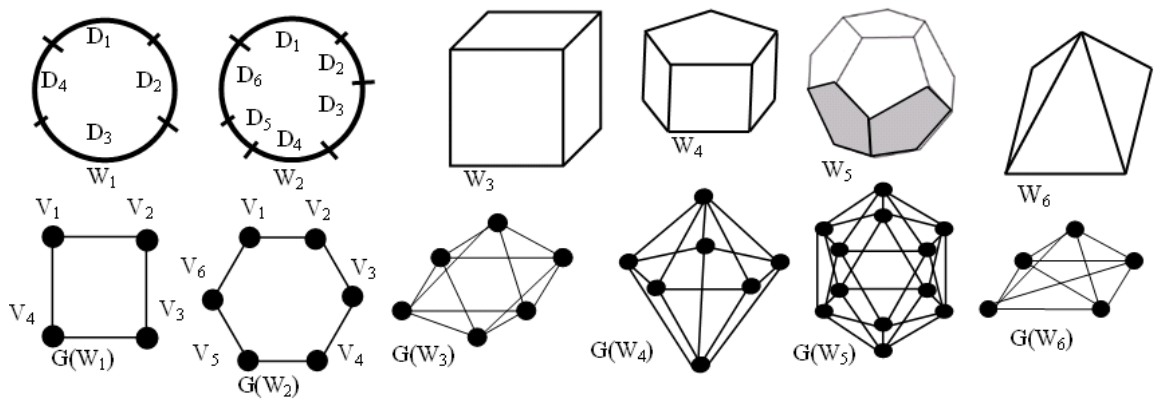

**Figure 11** $W_1$ and $W_2$ are LCL covers of a circle. Their intersection graphs $G(W_1)$ and $G(W_2)$ are digital 1-spheres. $W_3$, $W_4$ and $W_5$ are LCL covers of a sphere. Their intersection graphs $G(W_3)$, $G(W_4)$ and $G(W_5)$ are digital 2-spheres. $W_6$ is a non-LCL covers of a sphere. Its intersection graph $G(W_6)$ is not a digital 2-sphere.

**Theorem 4.2**

Let X be an n-dimensional sphere and an LCL collection $W=\{D_0,D_1...D_t\}$ of n-cells be a cover of $X=D_0 \cup D_1 \cup ...D_t$. Then the intersection graph $G(W)$ of W is a digital n-sphere.

Proof.
The proof is by induction on dimension n. For n=1, the statement is verified directly (see figure 11). Assume that the statement is valid whenever n<m. Let n=m. Let $V=\{D_1,D_2,...D_s\}$ be a collection of all n-cells intersecting $D_0$ (excluding $D_0$). Then the LCL collection of (n-1)-cells $V_1=\{C_1,C_2,...C_s\}$, where $C_i=D_0 \cap D_i$, i=1,2,...s, is a cover of the (n-1)-sphere $\partial D_0$, and collections V and $V_1$ are isomorphic according to theorem 2.1. Since the intersection graph $G(V_1)$ of $V_1$ is a digital (n-1)- sphere by the induction hypothesis, then G(V) is the digital (n-1)-sphere. For the same reason, the rim of every point in G(W) is a digital (n-1)-sphere. Therefore, G(W) is a digital n-manifold according to definition 3.2. Since $D_0$ is an n-cell then $E=X-intD_0$ is also an n-cell by fact 2.1. $U=\{D_1,D_2,…D_t\}$ is an LCL cover of $E= D_1\cup D_2\cup…D_t$. Therefore, the intersection graph $G(U)=\{v_1,v_2…v_t\}=G(W)-\{v_0\}$ is contractable by theorem 4.1. Hence, G(W) is a digital n-sphere by definition 3.1. This completes the proof. □

It is easy to check directly that the minimal LCL cover of an n-dimensional sphere contains 2n+2 n-cells, the minimal LCL cover of a torus and a Klein bottle contains sixteen 2-cells, the minimal LCL cover of a projective plane contains eleven 2-cells and the minimal LCL cover of a Moebius band contains twelve 2-cells.

**Example 4.1 Digital model of a circle**

Figure 11 illustrates theorem 4.2 for a circle. $W_1$ is an LCL cover of C. The intersection graph $G(W_1)$ is a digital 1-sphere. Notice that $W_1$ is an irreducible collection of 1-cells, and $G(W_1)$ is a digital minimal 1-sphere. $W_2$ is also an LCL cover of C, and $G(W_2)$ is a digital 1-sphere.

**Example 4.2 Digital model of a sphere**

Figure 11 depicts covers of a sphere. Cover $W_3$ consists of six square faces of a cube. $W_3$ is an irreducible LCL cover of a sphere. The intersection graph $G(W_3)$ is a minimal digital 2-sphere. Collection $W_4$ is composed of seven faces of a pentagonal prism. $W_4$ is an LCL cover of a sphere, $G(W_4)$ is a digital 2-sphere. $W_5$ is a set of twelve faces of a dodecahedron. Clearly, $W_5$ is an LCN cover of a sphere. Its intersection graph $G(W_5)$ is a digital 2-sphere. Faces of a square pyramid $W_6$ are a non-LCL cover of a sphere. Therefore, the intersection graph $G(W_6)$ is not a digital 2-sphere.

**Remark 4.1**

An irreducible LCL cover $W=\{D_1,D_2…D_{2n+2}\}$ of an n-sphere S by n-cells is the set of n-dimensional faces contained in the boundary $\partial U^{n+1}$ of an (n+1)-dimensional cube $U^{n+1}$ in the Euclidean (n+1)-dimensional space $E^{n+1}$. The boundary $\partial U^{n+1}$ topologically is an n-dimensional sphere. Obviously, n-dimensional faces $D_k$, k=1,2,…2n+2, of $U^{n+1}$ compose an LCL cover. The intersection graph G(W) is a digital minimal n-sphere containing 2n+2 points. It follows from examples 4.1 and 4.2 that a minimal digital model of an n-dimensional sphere contains 2n+2 points.

The proof of the following result is similar to the proof of theorem 4.2.

**Theorem 4.3**

>   Let X be a closed n-manifold and an LCL collection $W=\{D_0,D_1…D_t\}$ of n-cells be a cover of $X=D_0\cup D_1\cup…D_t$. Then the intersection graph G(W) of W is a digital n-manifold.

Proof.
Let $V=\{D_1,D_2,…D_s\}$ be a collection of all n-cells intersecting $D_0$ (excluding $D_0$ itself). It follows from the proof of theorem 4.2 that G(V) is a digital (n-1)-sphere. For the same reason, the rim of any point in G(W) is the digital (n-1)-sphere. Therefore, G(W) is a digital n-manifold by definition 2.1. This completes the proof. □

**Example 4.3 Digital model of a torus T, a projective plane P, Klein bottle K and a Moebius band M**

LCL covers of a torus T, a projective plane P, Klein bottle K and a Moebius band M and their intersection graphs are shown in figure 12. Covers W(T), W(P), W(K) and W(M) are irreducible. G(T), G(P), G(K) and G(M) are digital 2-manifolds with the minimal number of points. It can be checked directly that the Euler characteristic and the homology groups of T, P, K and M are the same as of G(T), G(P), G(K) and G(M) respectively.

One can ask what is the difference between two digital models of the same continuous closed n-dimensional manifold. In all examples, if $W_1$ and $W_2$ are different LCL covers of the same closed n-dimensional manifold M then their intersection graphs are homotopy equivalent, i.e., one of them can be converted into the other by a sequence of contractible transformations. To prove this fact is a problem for further study.

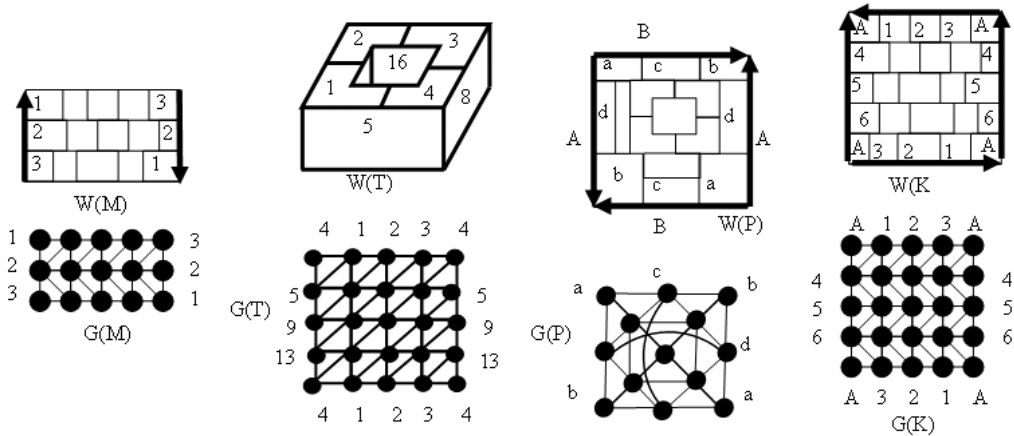

**Figure 12**. W(T), W(P), W(K) and W(M) are LCL covers of a torus T and a projective plane P, a Klein bottle K and a Moebius band M, and G(T), G(P), G(K) and G(M) are their intersection graphs. G(T), G(P), G(K) and G(M) are digital 2-manifolds.

## Remark 4.2

Next consider the following procedure for constructing a digital model of a continuous n-dimensional object in general and a continuous closed n-dimensional manifold in particular.
- Let X be a continuous closed n-dimensional manifold. Choose an LCL cover $W=\{D_1, D_2 \ldots D_s\}$ of X by n-tiles. Notice that the size and the form of any n-tile can be arbitrary depending on specific requirements but W must be an LCL collection.
- Build the intersection graph $G(W) = \{v_1, v_2 \ldots v_s\}$ of W. G(W) is necessarily a digital n-dimensional manifold with the same structural properties as X.

## Conclusion

- This paper studies properties of an LCL collection of n-cells on $E^p$, $p \geq n$. It is shown that a digital model of a continuous closed n-dimensional manifold, i.e., the intersection graph of an LCL cover of a continuous closed n-dimensional manifold, is necessarily a digital n-dimensional manifold; and that a digital model of a continuous n-dimensional sphere is necessarily a digital n-dimensional sphere containing at least 2n+2 points.
- An obvious consequence of this approach is that it demonstrates that in physical investigations and theories involving continuous n-dimensional manifolds such as n-dimensional spheres, a Möbius strip, a Klein bottle and so on, these manifolds can be replaced by their digital models, which are graphs with the same mathematical structure as their continuous counterparts. This essentially simplifies the study of processes and theories involving n-dimensional manifolds because this makes it possible to consider discrete spaces with the finite or countable number of points and to use a wide range of mathematical methods including computational and numerical methods in an obvious way.
- A problem for further study. In all examples, the intersection graphs of different LCL covers of a given continuous closed n-dimensional manifold are homotopy equivalent to each other, i.e., structural properties of a digital model do not depend on the choice of an LCL cover. Future research could be conducted to prove this fact.